\documentclass[12pt]{article}
\usepackage{amsmath}
\usepackage{amssymb}
\usepackage{latexsym}

\def\la {\lambda}
\def\coeff {{\Big|}}
\def\snake {\rfloor}
\def\Tau {{\mathcal T}}
\def\H {{\mathcal H}}
\def\V {{\mathcal V}}

\def\endofproof {\hskip .1in $\Box$}

\newcommand{\ontop}[2]{\substack{#1\\#2}}
\def\today{\ifcase\month\or
January\or February\or March\or April\or may\or June\or
July\or August\or September\or October\or November\or
December\fi
\space\number\day, \number\year}
\numberwithin{equation}{section}
\newtheorem{thm}{Theorem}[section]
\newtheorem{lemma}[thm]{Lemma}
\newtheorem{prop}[thm]{Proposition}

\newtheorem{cor}[thm]{Corollary}

\begin{document}

\title{Vertex operators for standard bases of the symmetric functions}
\author{Mike Zabrocki\\ Centre de Recherche Math\'ematiques, Universit\'e \\ de
Montr\'eal/LaCIM, Universit\'e du Qu\'ebec \`a Montr\'eal \\ email: {\small\tt
zabrocki@math.ucsd.edu}}
\date{}
\maketitle
\smallskip
\centerline{}
\smallskip
{\small MR Subject Number: {05E05}\hfill}

\smallskip
{\small Keywords: symmetric functions, vertex operators\hfill}
\bigskip

\begin{abstract}
We present a formulas to add a row or a column to the power, monomial,
forgotten, Schur, homogeneous and elementary symmetric functions.  As an application
of these operators we show that the operator that adds a column to the Schur
functions can be used to calculate a formula for the number of pairs of
standard tableaux the same shape and height less than or equal to a fixed $k$.
\end{abstract}
\vskip .3in

\section{Notation}

Using the notation of \cite{M},
we will consider the power
$\left\{ p_\la[X]
\right\}_\la$, Schur $\left\{ s_\la[X]
\right\}_\la$, monomial $\left\{ m_\la[X] \right\}_\la$, homogeneous $\left\{ h_\la[X]
\right\}_\la$, elementary $\left\{ e_\la[X] \right\}_\la$ and forgotten $\left\{ f_\la[X]
\right\}_\la$  bases for the symmetric functions.  We will often appeal to \cite{M} 
for proofs of identities relating these bases.

These bases are all indexed by partitions, non-increasing sequences of non-negative
integers.  The $i^{th}$ entry of the partition will be denoted by $\la_i$.
The length of a partition $\la$ is the largest $i$ such that $\la_i$ is non-zero
and will be denoted by $l(\la)$.
The size of the partition will be denoted by $|\la|$ and is equal to the sum over all
the entries of $\la$.  The symbol $n_i(\la)$ will be used to represent the number
of parts of size $i$ in the partition $\la$.  The conjugate partition will be
denoted by $\la'$ and is the partition such that $\la_i' =$ the number of $j$ such that
$\la_j$ is greater than or equal to $i$.

There is a standard inner product on symmetric functions $\left< p_\la, p_\mu \right> =
z_\la \delta_{\la\mu}$ (where $\delta_{xy} = 1$ if $x=y$ and $0$ otherwise and $z_\la =
\prod_{i \geq 1} i^{n_i(\la)} n_i(\la)!$.

We will use a few non-standard operations on partitions that will require
some notation.  The first is adding a column (or a sequence of columns) to
a partition.  Let $\la + a^k$ denote the partition $(\la_1 +a, \la_2 +a, \ldots,
\la_k + a)$.  We will assume that this partition is undefined when $l(\la)>k$.

Use the notation $\la - (\mu)$ to denote the partition formed by removing
the parts that are equal to $\mu$ from the partition $\la$.  This of course assumes
that there is a sequence $I = \{ i_1, i_2, \ldots i_{l(\mu)} \}
\subset \left\{1,2,\ldots l(\la) \right\}$ such that
$\mu = (\la_{i_1}, \la_{i_2}, \ldots, \la_{i_{l(\mu)}})$.  If this sequence does not exist
then $\la - (\mu)$ is again undefined.

The last operation will be inserting
parts into a partition and will be represented by $\la + (\mu)$.  This will be
the partition formed by ordering the sequence $(\la_1, \la_2, \ldots, \la_{l(\la)},
\mu_1, \mu_2, \\ \ldots, \mu_{l(\mu)})$ into a partition.

We will say that two bases for the symmetric functions
$\{a_\la\}_\la$ and $\{b_\la\}_\la$ are
dual if they have the property that $\left< a_\la, b_\mu \right> = \delta_{\la\mu}$.
By definition, the power symmetric functions are dual to the basis
$\{ p_\la/z_\la \}_\la$.  The monomial and homogeneous
symmetric functions are dual. The forgotten and the
elementary symmetric functions are dual.  The Schur symmetric functions
are self dual ( $\left< s_\la, s_\mu \right> = \delta_{\la\mu}$).

There exists an involution, $\omega$,
on symmetric functions that relates the elementary and
homogeneous bases by $\omega h_\mu = e_\mu$ and the monomial and forgotten bases
by $\omega m_\mu = f_\mu$.  It also has the property that $\omega s_\la = s_{\la'}$.

Denote the operation of 'skewing' by a symmetric function $f$ by $f^{\perp}$.
It is defined as being the operation dual to multiplication by the symmetric
function $f$ in the sense that $\left< f^\perp P, Q \right> = \left< P, f Q \right>$.
Its action on an arbitrary symmetric function $P$ may be calculated by the formula
$f^{\perp} P = \sum_\la \left< P, f a_\la \right> b_\la$ for any dual bases $\{a_\la\}_\la$
and $\{b_\la\}_\la$.

Using results and notation in \cite{M} (p.92-3 example (I.5.25)),
by setting $\Delta f = \sum_\mu (a_\mu^\perp f) \otimes b_\mu^\perp$
where $\{ a_\mu \}_\mu$ and $\{ b_\mu \}_\mu$ are any dual bases.  It follows that
if $\Delta f = \sum_i c_i \otimes d_i$ then $f^\perp(PQ) = \sum_i c_i^\perp(P) d_i^\perp(Q)$.
Using this definition, $\Delta h_k = \sum_i h_i \otimes h_{k-i}$ and
$\Delta e_k = \sum_i e_i \otimes e_{k-i}$ and $\Delta p_k = p_k \otimes 1 + 1 \otimes p_k$.

By the phrase 'vertex operators' we are referring to linear symmetric function operators that
add a row or a column to the partitions indexing a particular family of symmetric functions. 
Formulas of this type for symmetric functions are sometimes called Rodrigues formulas.  In this
article we look at those symmetric function operators which lie in the linear span of 
$\left\{ f_i g_i^\perp \right\}_i$ where $f_i$ and $g_i$ are symmetric functions
to find expressions for vertex operators for each basis.  

The vertex operators for the elementary
and forgotten symmetric function basis are related to the operators for the homogeneous
and monomial (resp.) symmetric functions by conjugating by the operator
$\omega$.

The existence of such operators for the Schur (\cite{M} p.96-7, \cite{Ze} p.69),
and (row only)
Hall-Littlewood (\cite{M} p.237-8, \cite{J}) symmetric functions are known.
For the multiplicative
bases, it is clear that there exists operators that add a row to the symmetric
functions of this form since $e_k e_\la = e_{\la+(k)}$, $h_k h_\la = h_{\la+(k)}$
and $p_k p_\la = p_{\la+(k)}$, but adding a
column is not an obvious operation.

In general, formulas for adding a row or a column can be useful in proving a
combinatorial interpretation for a symmetric function or deriving new formulas
or properties.  The author's interest in this particular question comes from trying to find
vertex operators for the Macdonald symmetric functions.  The Macdonald 
vertex operator must specialize to the vertex operators for other symmetric functions
and so understanding these operators is an important first step.

\section{The power vertex operator}

This is the warm up case for the other $5$ bases.  The commutation
relation between
$p_k^\perp$ and $p_j$ has a nice expression: $p_k^\perp p_j = p_j p_k^\perp
+ k \delta_{kj}$.  This can be used to show the slightly more general relation $p_\la^\perp p_k =
p_k p_\la^\perp + k n_k(\la) p_{\la-(k)}^\perp$ (where it is assumed that $p_{\la-(k)}^\perp = 0$
if $\la$ does not contain a part of size $k$).

The vertex operator is given by the following theorem

\begin{thm} For $k\geq 0$ define the following linear operator
$$CP_{a^k} = \sum_{\ontop{\la}{l(\la) \leq k}} p_a^{k-l(\la)} \prod_{i=1}^{l(\la)}
\left( p_{\la_i+a} - p_{\la_i} p_a \right) p_\la^\perp/z_\la$$
where the sum is over all partitions $\la$ with less than or equal to $k$ parts (if
$k=0$ then $CP_{a^0}=1$).
$CP_{a^k}$ has the property that $CP_{a^k} p_\mu = p_{\mu+a^k}$ for all $\mu$ such
that $l(\mu)<k$.
\end{thm}
\vskip .3in

\noindent
{\bf Proof:}
The proof is by induction on the number of parts of $\mu$.
Clearly this operator has the property that $CP_{a^k} 1 = p_a^k$ since $p_\la^\perp$
kills $1$ for $|\la|>0$.  From the
commutation relation of $p_\la^\perp$ and $p_k$ we derive that
$$CP_{a^k} p_j = (p_{j+a} - p_j p_a) CP_{a^{k-1}} + p_{j} CP_{a^k}$$
The proof by induction follows from this relation.
\endofproof
\vskip .3in

The formula $CP_{a^k}$ was chosen so that it has two properties: it adds a column to the
power symmetric functions, and it has a relatively simple expression when written in this
notation.  The action of this operator on $p_\la$ when $l(\la)>k$ is not specified by
these conditions, but it is determined.

If one wishes to give an expression for an operator that has the same
action on $p_\la$ for $l(\la) \leq k$ and the action on $p_\la$ for $l(\la)>k$ is something
else (say for instance $0$), this is possible by adding in terms of the form $p_\mu p_\la^\perp$
where $l(\la)>k$ to $CP_{a^k}$.

\section{Homogeneous and elementary vertex operators}

Note that 
$\left< h_i^\perp m_\la, h_\mu \right>
= \left< m_\la, h_i h_\mu \right> = \delta_{\mu, \la - (i)}$.  Therefore 
\begin{equation}
h_i^\perp (m_\la) = m_{\la -(i)} \label{hponm}
\end{equation} and $h_i^\perp (m_\la) = 0$ if $\la$ does not have a part of size $i$.

We then use these results to prove the following lemma
\begin{lemma} \label{hkmlacomm}
$$h_k^\perp m_\la = \sum_{i \geq 0}
m_{\la - (i)} h_{k-i}^\perp$$
$$m_\la^\perp h_k = \sum_{i \geq 0}
h_{k-i} m_{\la - (i)}^\perp$$
where we will assume the convention $m_{\la -(i)} = 0$ whenever $\la -(i)$ is undefined.
\end{lemma}
\vskip .3in

\noindent
{\bf Proof:}

The first identity follows from the remarks made in the previous section
that say that
$h_k^\perp (m_\la P) = \sum_i h_i^\perp (m_\la) h_{k-i}^\perp (P)$ (where $P$ is any
arbitrary symmetric function). 

The second identity is a restatement of the first since
\begin{equation*}
\left< m_\la^\perp h_k P, Q \right> = \left< P, h_k^\perp m_\la Q \right>
= \sum_{i \geq 0} \left< P, m_{\la - (i)} h_{k-i}^\perp Q \right>  
= \sum_{i \geq 0} 
\left< h_{k-i} m_{\la - (i)}^\perp P, Q \right>
\end{equation*}
\endofproof
\vskip .3in

Define $CH_{1^k}$ to be the operator $CH_{1^k} = \sum_\la (-1)^{|\la|} e_{\la+1^k} m_\la^{\perp}$,
and the operator $CE_{1^k}$ to be $CE_{1^k} = \sum_\la (-1)^{|\la|} h_{\la+1^k} f_\la^{\perp}$
where we will assume that $h_{\la+1^k} = 0$ if $l(\la)>k$ so that the sums
in these equations are over all partitions with parts 
smaller than or equal to $k$.

\vskip .2in
The vertex operator property that we prove for the homogeneous and elementary symmetric
functions is
\begin{thm}
If $l(\la) \leq k$, then $CH_{1^k} h_\la = h_{\la+1^k}$ and $CE_{1^k} e_\la = e_{\la+1^k}$ .
\end{thm}
\vskip .3in

\noindent
{\bf Proof:} 
The proof is a matter of showing that for $k>0$ operator $CH_{1^k}$ and $h_n$ 
(considered as an operator that consists of multiplication by $h_n$) has
the commutation relation $CH_{1^k} h_n = h_{n+1} CH_{1^{k-1}}$.
$$CH_{1^k} h_n = \sum_\la (-1)^{|\la|} e_{\la+1^k} m_\la^{\perp} h_n$$
The sum here is over $\la$ with the number of parts less than or equal to $k$.
Apply Lemma \ref{hkmlacomm} and rearrange the terms in the sum.
\begin{align*}
CH_{1^k} h_n &= \sum_\la (-1)^{|\la|} e_{\la+1^k} \sum_{i \geq 0}
h_{n-i} m_{\la - (i)}^\perp \\
&= \sum_\la \sum_{i \geq 0}
(-1)^{|\la|-i} (-1)^i e_{\la+1^k} h_{n-i} m_{\la - (i)}^\perp\\
&= \sum_\la \sum_{i \geq 0}
(-1)^{|\la|-i} (-1)^i h_{n-i} e_{i+1} e_{\la- (i) +1^{k-1}} m_{\la - (i)}^\perp
\end{align*}

In the last equation, there is an assumption that $e_{\la- (i) +1^{k-1}}=0$
if $\la - (i)$ is undefined.
As long as $k>1$, making the substitution $\la \rightarrow \la+(i)$ yields the
equation: 
$$=\left( \sum_{i =1}^n
(-1)^i h_{n-i} e_{i+1} \right) \left( \sum_\la (-1)^{|\la|} 
e_{\la +1^{k-1}} m_{\la}^\perp \right)$$
Now the sum is over $\la$ with less than or equal to $k-1$ parts.  This is equal to
$=h_{n+1} CH_{1^{k-1}}$
using the well known relation $\sum_{r=0}^n (-1)^r e_r h_{n-r} = 0$ for $n\geq 0$. 
If $k=1$ then 
\begin{align*}
CH_{1^1} h_n &= \sum_{i =1}^n (-1)^i h_{n-i} e_{i+1} + \sum_{i \geq 1} (-1)^{i} e_{i+1} h_n
m_{(i)}^{\perp} \\ &= h_{n+1} + \sum_{i \geq 1} (-1)^{i} e_{i+1} h_n m_{(i)}^{\perp}
\end{align*}
which is the 'correct' answer only if $CH_{1^1} h_n$ are acting on $1$.  Notice also that
$CH_{1^k}$ acting on $1$, yields $h_1^k$ since only one term is not $0$.

The corresponding result for the $CE_{1^k}$ operator follows by noting that
$CE_{1^k} = \omega CH_{1^k} \omega$.
\endofproof\vskip .3in

The action of $CH_{1^k}$ on $h_\la$ when $l(\la) > k$ is not known.  
The sum in the formula for $CH_{1^k}$ is only over partitions $\la$ such
that $l(\la) \leq k$ and by adding terms of the same form but with $l(\la)>k$
it is possible to modify the formula so that the action on the
$h_\la$ when $l(\la)>k$ is $0$, but the formula will not be as simple.

It would be interesting to know the action of these vertex operators on other bases
besides the one that it adds a row and column to.  For instance, actions of $e_k$, $h_k$, and
$p_k$ are known on the Schur basis, but what is the action
of an operator that adds a column to the homogeneous, elementary, or power basis when
it acts on the Schur basis?

Note the following two formulas that relate $CH_{1^k}$ and $CE_{1^k}$.

\begin{equation} \label{HErel}
CH_{1^k} = \sum_{\ontop{\la}{l(\la) \leq k}} (-1)^{|\la|} CE_{1^k} (e_\la) m_\la^\perp
\end{equation}
\begin{equation} \label{EHrel}
CE_{1^k} = \sum_{\ontop{\la}{l(\la) \leq k}} (-1)^{|\la|} CH_{1^k} (h_\la) f_\la^\perp
\end{equation}

This is the first instance when a pair of operators share a relation like this,
and it will occur with pairs of the other operators that appear in this article.
These relations fall under the category of 'eerie coincidences' (by this I mean that
there is probably some explanation for these relations but they are very unexpected
and I don't know what that explanation might be).

\section{Monomial and forgotten vertex operators}

The vertex operators for the monomial and forgotten symmetric functions requires
a few identities.

\begin{lemma} \label{expan}
Let $r_{\mu} = (-1)^{|\mu| - l(\mu)} \frac{l(\mu)!}{n_1(\mu)! n_2(\mu)! \cdots}$
then for $|\mu| \geq 0$,
$e_k = \sum_{\mu \vdash k} r_{\mu} h_\mu$
\end{lemma}

\noindent
{\bf Proof: }
\cite{M} example I.2.20 p.33 \endofproof
\vskip .3in

\begin{lemma} \label{ecid}
$\sum_{j \geq 0} (-1)^j r_{\mu - (j)} = 0$ where is is assumed that
if $\mu-(j)$ does not exist then $r_{\mu-(j)}=0$.
\end{lemma}

\noindent
{\bf Proof:}
$\sum_{j=0}^k (-1)^j e_{k-j} h_j = 0$.  Now expand $e_{k-j}$ in terms of the
homogeneous basis using the last lemma and equate coefficients of $h_\mu$
on both sides of the equation.
 \endofproof
\vskip .3in

\begin{lemma} \label{ekmlacomm}
$e_k^\perp m_\la = \sum_\mu r_{\mu} m_{\la - (\mu)} e_{k-|\mu|}^\perp$
where $m_{\la - (\mu)} = 0$ if $\la - (\mu)$ is undefined.
\end{lemma}

\noindent
{\bf Proof:}
By \cite{M} (example I.5.25, p.92-3), $e_k^\perp m_\la =
\sum_{i \geq 0} e_i^\perp (m_\la) e_{k-i}^\perp$.
The expansion of the $e_i^\perp$ in terms of $h_\mu^\perp$ is given in the last lemma and so
we have that 
\begin{equation*}
e_k^\perp m_\la = \sum_{i \geq 0} \left(\sum_{\mu \vdash i} r_{\mu} h_\mu^\perp (m_\la)
\right) e_{k-i}^\perp
= \sum_{i \geq 0} \left(\sum_{\mu \vdash i} r_{\mu} m_{\la-(\mu)}
\right) e_{k-i}^\perp
\end{equation*}
by the equation (\ref{hponm}) and this is equivalent to the statement of the lemma.\endofproof
\vskip .3in

\begin{lemma} \label{mkmla}
$m_{(k)} m_\la = \sum_{i \geq 0} (1+n_{k+i}(\la)) m_{\la -(i) + (k+i)}$
where it is assumed that $m_{\la -(i) + (k+i)} = 0$ if
$\la - (i)$ is undefined.
\end{lemma}

\noindent
{\bf Proof:}  For paritions $\mu$ of $|\la|+k$, one has that 
$$m_{(k)} m_\la \coeff_{m_\mu} = h_\mu^\perp (m_{(k)} m_\la)$$

We note that for all $n \geq 0$ that
$h_n^\perp m_{(k)} = m_{(k)} h_n^\perp + h_{n-k}^\perp$.  Apply this to
the expression for the coefficient of $m_{\mu}$
$$h_\mu^\perp (m_{(k)} m_\la) = \sum_{j=1}^{l(\mu)} h_{\mu - (\mu_j) + (\mu_j-k)}^\perp( \mu_\la)$$

This implies that for the coeffient to be non-zero that $\mu$ must be equal to $\la$ with
a part (say of size $i$) pulled away and a part of size $k+i$ added in.  The coefficient
will be the number of times that $k+i$ appears in the partition $\mu$ (one more time
than it appears in the partition $\la$).
\endofproof
\vskip .3in

The first vertex operator that is presented here
for the monomial symmetric functions adds a
row but it also multiplies by a coefficient
(a property that is unwanted in the final result), but this operator
provides an easy method for obtaining an operator that does not have
this coefficient.
\begin{prop}
Let $RM_a^{(1)} = \sum_{i\geq 0} (-1)^i m_{(a+i)} e_{i}^\perp$
then 
$RM_a^{(1)} m_\la = (1 + n_a(\la)) m_{\la + (a)}$
\end{prop}
\vskip .3in

\noindent
{\bf Proof:} 
$$RM_a^{(1)} m_\la
= \sum_{i\geq 0} (-1)^i m_{(a+i)} e_{i}^\perp m_\la$$

Apply Lemma \ref{ekmlacomm} to get
$$= \sum_{i\geq 0} (-1)^i m_{(a+i)}
\sum_{\mu \vdash i} r_{\mu} m_{\la - (\mu)}
$$

The sum over $i$ and $\mu$ may be combined to form 
one sum over all partitions $\mu$.

$$= \sum_{\mu} (-1)^{|\mu|} r_{\mu}
m_{(a+|\mu|)} m_{\la - (\mu)}
$$

Now multiplying by a monomial symmetric function with
one part has an expansion given in Lemma \ref{mkmla}.

$$= \sum_{\mu} \sum_{j \geq 0}
(-1)^{|\mu|} r_{\mu}
(1+n_{a+|\mu|+j}(\la - (\mu))) m_{\la - (\mu) - (j) + (j+a+|\mu|)}
$$

The terms indexed by the same monomial symmetric function may be
grouped together by letting $\nu = \mu + (j)$.

$$=\sum_\nu \sum_{j\geq 0} (-1)^{|\nu| - j} r_{\nu-(j)}
(1+n_{a+|\nu|}(\la - ((\nu)-(j)))) m_{\la - (\nu) + (a+|\nu|)}$$

$$=\sum_\nu (1+n_{a+|\nu|}(\la)) m_{\la - (\nu) + (a+|\nu|)}
\sum_{j\geq 0} (-1)^{|\nu| - j} r_{\nu-(j)}$$

But $\sum_{j\geq 0} (-1)^{|\nu| - j} r_{\nu-(j)} = 0$
if $|\nu|>0$ by Lemma \ref{ecid}.  There is one term left.

$$=(1+n_{a}(\la)) m_{\la + (a)}$$
\endofproof
\vskip .3in

An expression for an operator that adds a row without a coefficient
can be written in terms of this operator.

\begin{thm}\label{Mavertex}
For $a>0$
define $RM_a = \sum_{k \geq 0} (-1)^k \frac{\left( RM_a^{(1)} \right)^{k+1} }{(k+1)!}
\left(h_a^k\right)^\perp$
then $RM_a m_\la = m_{\la+(a)}$.
\end{thm}
\vskip .3in

\noindent
{\bf Proof:}
Apply the previous proposition to this formula and reduce using the following steps.
\begin{align*}
RM_a m_\la &= \sum_{i \geq 0} (-1)^i
\frac{\left( RM_a^{(1)} \right)^{i+1} }{(i+1)!}
\left(h_a^i\right)^\perp m_\la \\
&= \sum_{i \geq 0} (-1)^i \frac{(n_{a}(\la)+1)(n_{a}(\la)+2) \cdots (n_{a}(\la)+i+1)}{(i+1)!}
m_{\la - (a^i) +(a^{i+1})} \\
&= \sum_{i = 0}^{n_{a}(\la)} (-1)^i \left( \begin{array}{c} n_{a}(\la)+1 \\ i+1 \end{array} \right)
m_{\la + (a)} \\
&= m_{\la + (a)}
\end{align*}

The last equality is true because for $a>0$,
$\sum_{i = 1}^{a} (-1)^{i-1} \left( \begin{array}{c} a \\ i \end{array} \right) = 1$.
\endofproof
\vskip .3in

Notice that the action of the $RM_a$ operators on the monomial basis implies
that $RM_a RM_b = RM_b RM_a$.  This property is difficult to derive just from the
definition of the operator.

This expression for the operator $RM_a$ is a little unsatisfying
since the computation of $\left( RM_a^{(1)} \right)^i$ can be simplified.  The
following operator shows how $RM_a$ can be reduced to closer resemble $CH_{1^k}$.
To add more than one row at a time to a monomial symmetric function, the formula resembles the
vertex operator that adds a column to the homogeneous basis.  

\begin{prop} \label{makop}
For $a \geq 0$ and $k \geq 0$, we have that
$$RM_a^{(k)} = \frac{\left( RM_a^{(1)} \right)^k}{k!} 
= \sum_\la (-1)^{|\la|} m_{\la + a^k} e_\la^\perp$$
with the understanding that $m_{\la + a^k} = 0$ if ${\la + a^k}$ is undefined and
$n_0(\la) = k - l(\la)$.
It follows that
$RM_a^{(k)} m_\la = \left( \begin{array}{c} n_a(\la) + k \\ k \end{array} \right) m_{\la + (a^k)}$.
\end{prop}
\vskip .3in

\noindent
{\bf Proof:}

By induction on $k$, we will show that $RM_a RM_a^{(k)} = (k+1) RM_a^{(k+1)}$.  It follows
that $(RM_a^{(1)})^k = k! RM_a^{(k)}$.  Since $(RM_a^{(1)})^k m_\la =
(n_a(\la) + 1) (n_a(\la) + 2) \cdots (n_a(\la) + k) m_{\la + (a^k)}$ then
$RM_a^{(k)} m_\la = \frac{(n_a(\la) + 1) (n_a(\la) + 2) \cdots (n_a(\la) + k)}{k!} m_{\la + (a^k)}$.

$$RM_a RM_{a}^{(k)} = 
\sum_{j \geq 0} (-1)^j m_{(a+j)} e_j^\perp \sum_\la (-1)^{|\la|} m_{ \la + a^k} e_\la^\perp$$

Commute the action of $e_j^\perp$ and $m_{ \la + a^k}$ using Lemma \ref{ekmlacomm}.

\begin{align*}
&=\sum_{j\geq 0} (-1)^j m_{(a+j)} \sum_\la (-1)^{|\la|} \sum_\mu r_{\mu} m_{\la + a^k - (\mu)}
e_{j-|\mu|}^\perp e_\la^\perp \\
&=\sum_{j \geq 0} \sum_\la \sum_\mu (-1)^{|\la|+j} r_{\mu} m_{(a+j)} m_{\la + a^k - (\mu)}
e_{j-|\mu|}^\perp e_\la^\perp
\end{align*}

The formula for multiplying a monomial symmetric function with one part is given in
Lemma \ref{mkmla}.

$$=\sum_{j\geq 0} \sum_\la \sum_\mu \sum_{l \geq 0}
(-1)^{|\la|+j} r_{\mu} n_{a+j+l}(\la + a^k - (\mu) - (l) + (a+l+j))
m_{\la + a^k - (\mu) - (l) + (a+l+j)} e_{j - |\mu|}^\perp e_\la^\perp$$

The next step is to change the sum over $\mu$ so that  it includes the part of
size $l$, this is equivalent to making the replacement $\mu \rightarrow \mu - (l)$.

$$=\sum_{j \geq 0} \sum_{\la} \sum_{\mu} \sum_{l \geq 0}
(-1)^{|\la|+j} r_{\mu-(l)} n_{a+j+l}(\la + a^k - (\mu) + (a+l+j))
m_{\la + a^k - (\mu) + (a+l+j)} e_{j+l - |\mu|}^\perp e_\la^\perp$$

Let $i=j+l$, then the sum over $j$ can be converted to a sum over $i$.

$$=\sum_{\la} \sum_{\mu} \sum_{l \geq 0} \sum_{i \geq l} 
(-1)^{|\la|+i-l} r_{\mu-(l)} n_{a+i}(\la + a^k - (\mu) + (a+i))
m_{\la + a^k - (\mu) + (a+i)} e_{i - |\mu|}^\perp e_\la^\perp$$

The next step is to interchange the sum over $i$ and the sum over $l$.  Since
$l \geq 0$ and $i \geq l$ then $i \geq 0$ and $0 \leq l \leq i$.

$$=\sum_{\la} \sum_{\mu} \sum_{i \geq 0} \sum_{l = 0}^{i} (-1)^l r_{\mu-(l)}
(-1)^{|\la|+i}  n_{a+i}(\la + a^k - (\mu) + (a+i))
m_{\la + a^k - (\mu) + (a+i)} e_{i - |\mu|}^\perp e_\la^\perp$$

Notice that since $e_{i - |\mu|}^\perp$ is zero for all $i<|\mu|$, then
all terms are zero unless $i \geq |\mu|$.

$$=\sum_{\la} \sum_{\mu} \sum_{i \geq |\mu|} \sum_{l = 0}^{i} (-1)^l r_{\mu-(l)}
(-1)^{|\la|+i}  n_{a+i}(\la + a^k - (\mu) + (a+i))
m_{\la + a^k - (\mu) + (a+i)} e_{i - |\mu|}^\perp e_\la^\perp$$

The sum over $l$ is equal to $0$ as long as $|\mu|>0$ using Lemma \ref{ecid}.
A substitution of $i \rightarrow i+|\mu|$ can be made so that $i$ runs over all
integers greater than or equal to $0$.

$$=\sum_{\la} \sum_{i \geq 0} (-1)^{|\la|+i}  n_{a+i}(\la + a^k + (a+i))
m_{\la + a^k + (a+i)} e_{i}^\perp e_\la^\perp$$

Let the sum over $\la$ include the part of size $i$, then $\la = \la + (i)$.

$$=\sum_{\la} (-1)^{|\la|} \sum_{i \geq 0} n_{a+i}(\la + a^{k+1})
m_{\la + a^{k+1}} e_\la^\perp$$

The sum over $i$ is now independent of $\la$ since $\sum_{i \geq 0} n_{a+i}(\la + a^{k+1})$
will always be $k+1$.

$$=(k+1) \sum_{\la} (-1)^{|\la|} m_{\la + a^{k+1}} e_\la^\perp = (k+1) RM_a^{(k+1)}$$\endofproof
\vskip .3in

It follows that the formula for $RM_a^{(k)}$ can be substituted into
Theorem \ref{Mavertex} and this provides a more reduced form of the
first formula given for $RM_a$.
\begin{cor}
$RM_a = \sum_{k \geq 0} \sum_\la (-1)^{|\la|+k} m_{\la + a^{k+1}} e_\la^\perp (h_a^k)^\perp$
\end{cor}
\vskip .3in

Since the forgotten basis is related to the monomial basis by an application of the
involution $\omega$, then the formulas for the symmetric function operator that
adds a row to the forgotten symmetric functions follows immediately.

\begin{cor}
$RF_a = \sum_{k \geq 0} \sum_\la (-1)^{|\la|+k} f_{\la + a^{k+1}} h_\la^\perp (e_a^k)^\perp$
has the property that
$$RF_a f_\la = f_{\la + (a)}$$
\end{cor}
\vskip .3in

There exists an operator $\Tau_{-X}$ of the same form as the operators
that exist already in this paper that has the property that $\Tau_{-X} P[X] =0$,
if $P[X]$ is a homogeneous symmetric function of degree greater than 0 and 
$\Tau_{-X} 1 = 1$.  This means that the operator applied to an arbitrary symmetric
function has the property that it picks out the constant term of the symmetric function.

\begin{prop}
Define the operator 
$$\Tau_{-X}= \sum_\la (-1)^{|\la|} s_\la s_{\la'}^\perp$$
Then for any dual bases $\{ a_\mu \}_\mu$ and $\{ b_\mu \}_\mu$ (that is,
$\left< a_\mu, b_\la \right> = \delta_{\la \mu}$), this is equivalent to
$$\Tau_{-X}= \sum_\la (-1)^{|\la|} \omega(a_\la) b_\la^\perp$$
This operator has the property that $\Tau_{-X} s_\la = 0$ for $|\la| > 0$
and $\Tau_{-X} 1 = 1$. 
\end{prop}
\vskip .3in

\noindent
{\bf Proof:}
$$\Tau_{-X} s_\mu = \sum_\la (-1)^{|\la|} s_\la s_{\la'}^\perp s_\mu$$
This is exactly the same expression as formula (\cite{M}, p.90, (I.5.23.1))
with the $x$ variables substituted for the $y$.  This expression is $0$
unless $s_\mu = 1$.

It requires very little to show that
this operator can be given an expression in terms of any dual
basis.
\begin{align*}
\Tau_{-X} &= \sum_\la (-1)^{|\la|} \omega(s_\la) s_{\la}^\perp \\
&= \sum_\la (-1)^{|\la|} \sum_{\mu\vdash |\la|} \left< \omega(s_\la), \omega(b_\mu) \right> \omega(a_\mu) s_{\la}^\perp \\
&=  \sum_\mu \sum_{\la \vdash |\mu|} (-1)^{|\mu|} \omega(a_\mu) \left< s_\la, b_\mu \right> s_{\la}^\perp \\
&= \sum_\mu (-1)^{|\mu|} \omega(a_\mu) b_\mu^\perp
\end{align*}\endofproof
\vskip .3in

Note that $\Tau_{-X}$ is actually a special case of a plethystic operator
$\Tau_{Z} P[X] = P[X+Z]$.

Fix a basis of the symmetric functions, $\{ a_\mu \}_\mu$, we may talk about the
symmetric function linear operator that sends $a_\mu$ to the expression $d_\mu$
(where $\{ d_\mu \}_\mu$ is any family of symmetric function expressions).
We can say that this operator lies in the linear span of the operators
$s_\la s_\mu^\perp$ and an expression can be given fairly easily.

\begin{cor} (The everything operator)

Let $\{ a_\mu \}_\mu$ be a basis of the symmetric functions and $\{ b_\mu \}_\mu$ be its dual
basis.  Then an operator that sends $a_\mu$ to the expression $d_\mu$ is
given by
$$E_{\{ a_\mu \}}^{\{ d_\mu \}} = \sum_\mu d_\mu \Tau_{-X} b_\mu^\perp$$
In other words we have that $E_{\{ a_\mu \}}^{\{ d_\mu \}}$ acts linearly,
and on the basis $a_\mu$ it has the action $E_{\{ a_\mu \}}^{\{ d_\mu \}} a_\mu = d_\mu$.
\end{cor}
\vskip .3in

\noindent
{\bf Proof:}

Note that when $b_\mu^\perp$ acts on a homogeneous polynomial, the result
is a homogeneous polynomial of degree $|\mu|$ less.  Therefore
if $|\mu| > |\la|$, then $b_\mu^\perp a_\la = 0$.  If $|\mu| < |\la|$ then
$\Tau_{-X} b_\mu^\perp a_\la = 0$ since $\Tau_{-X}$ kills all non-constant
terms.  When $|\mu| = |\la|$, we have that $b_\mu^\perp a_\mu = \delta_{\la\mu}$
and therefore, $\Tau_{-X} b_\mu^\perp a_\la  = \delta_{\la\mu}$.  This
also implies that $$\sum_\mu d_\mu \Tau_{-X} b_\mu^\perp a_\la = \sum_\mu d_\mu \delta_{\la\mu} =
d_\la$$\endofproof
\vskip .3in

This operator looks too general to be of much use, but using known
symmetric function identities it is possible to reduce and derive
expressions for other operators.  For instance, the
symmetric function operator that adds a column to the monomial
symmetric functions is a special case of this.

\begin{thm} \label{colum}

Let $CM_{a^k} = \sum_\la (-1)^{|\la|} \left( \begin{array}{c} n_a(\la) + k \\
k \end{array} \right) m_{\la + (a^k)} e_\la^\perp$, then
$CM_{a^k} m_\la = m_{\la + a^k}$ with the convention that
$m_{\la + a^k} = 0$ if $\la + a^k$ is undefined.
\end{thm}
\vskip .3in

\noindent
{\bf Proof:}

We will reduce an expression for $E_{\{ m_\la \}}^{\{ m_{\la + a^k} \}}$ to one
for $CM_{a^k}$.

$$E_{\{ m_\la \}}^{\{ m_{\la + a^k} \}} =
\sum_\la m_{\la + a^k} \sum_\mu (-1)^{|\mu|} m_\mu e_\mu^\perp h_\la^\perp$$

Let $r_{\la \mu}$ be the coefficient of $e_\mu$ in $h_\la$ (by an application of
the involution $\omega$ it is also the coefficient of $h_\mu$ in $e_\la$).
Then the expression is equivalent to

$$=\sum_\la m_{\la + a^k} \sum_\mu (-1)^{|\mu|} m_\mu 
e_\mu^\perp \sum_{\gamma \vdash |\la|} r_{\la \gamma} e_\gamma^\perp$$

Rearranging the sums this may be rewritten as
$$=\sum_\la \sum_\mu \sum_{\gamma} (-1)^{|\mu|} m_{\la + a^k}
r_{\la \gamma} m_\mu e_\mu^\perp e_\gamma^\perp$$

It is possible to group all the terms that skew by the same elementary symmetric function
by making the substitution $\mu \rightarrow \mu - (\gamma)$ since the sum over
$\mu$ and $\gamma$ are over partitions.

$$=\sum_\la \sum_\mu \sum_{\gamma} (-1)^{|\mu|-|\gamma|} m_{\la + a^k}
r_{\la \gamma} m_{\mu-(\gamma)} e_\mu^\perp$$

Note that $m_{\mu-(\gamma)} = h_\gamma^\perp (m_\mu)$ and $\sum_\gamma (-1)^{|\gamma|} r_{\la
\gamma} h_\gamma^\perp = (-1)^{|\la|} e_\la^\perp$.

$$=\sum_\mu \sum_\la (-1)^{|\mu|} m_{\la + a^k}
(-1)^{|\la|} e_\la^\perp (m_{\mu}) e_\mu^\perp$$

Notice that the first part of this expression is exactly the operator $RM_a^{(k)}$ acting
exclusively on $m_\mu$.  We may then apply Proposition \ref{makop} and note that the
expression reduces to the sum stated in the theorem.
\endofproof
\vskip .3in

The symmetric function operator that adds a column (or a group of columns) to
the forgotten symmetric functions can be found by conjugating the $CM_{a^k}$ operator
by the involution $\omega$ to derive the following corollary.

\begin{cor}

Let $CF_{a^k} = \sum_\la (-1)^{|\la|} \left( \begin{array}{c} n_a(\la) + k \\
k \end{array} \right) f_{\la + (a^k)} h_\la^\perp$, then
$CF_{a^k} f_\la = f_{\la + a^k}$ with the convention that
$f_{\la + a^k} = 0$ if $\la + a^k$ is undefined.
\end{cor}

The operator that adds a sequence of rows to the monomial symmetric functions and
the operator that adds a sequence of columns are related by a pair of formulas
similar to in the case of formulas (\ref{HErel}) and (\ref{EHrel}).  Notice
that proposition \ref{makop} and theorem \ref{colum} say that

\begin{equation} \label{CMrel}
CM_{a^k} = \sum_\la (-1)^{|\la|} RM_a^{(k)}( m_{\la}) e_\la^\perp
\end{equation}
\begin{equation} \label{MCrel}
RM_a^{(k)} = \sum_\la (-1)^{|\la|} CM_{a^k}(m_{\la}) e_\la^\perp
\end{equation}

This is 'eerie coincidence' number two.  The relation between these two operators is
very similar to the relation between $CH_{1^k}$ and $CE_{1^k}$ but not exactly the same.
Once again this is unexpected and unexplained.

\section{Schur vertex operators}

A symmetric function operator that adds a row to the Schur functions is given in \cite{M}
(p.95-6 I.5.29.d) that is of the same flavor as the other vertex operators presented here.

\begin{thm} (Bernstein) \label{bern}

Let $RS_a = \sum_{i\geq 0} (-1)^i h_{a+i} e_i^\perp$, then $RS_a s_\la = s_{\la+(a)}$
if $a \geq \la_1$.  In addition, $RS_a RS_b = - RS_{b-1} RS_{a+1}$.
\end{thm}

\noindent {\bf Proof:}
Repeated applications of this operator yeilds
expressions of the Jacobi-Trudi sort.
Use the relation $RS_a h_k = h_k RS_a - h_{k-1} RS_{a+1}$ (which follows
from \cite{M} example (I.5.29.b.5) and (I.5.29.d)),
$RS_a( 1 )= h_a$ and follow the proof of \cite{M} (I.3.(3.4'') p.43) which does not
actually require that the indexing sequence be a partition.
It follows then that
$$RS_{s_1} RS_{s_2} \cdots RS_{s_n} (1) = det \left| h_{s_j-j+i} \right|_{1\leq i,j \leq n}$$
\endofproof
\vskip .3in

Conjugating this operator by $\omega$ produces an operator that adds
a column to a Schur symmetric function.  We will show in this section that
an nice expression exists for a formula for an operator that adds a column
to a Schur function, but with the property that the result is $0$ if
the partition is longer than the column being added.

It follows from the commutation relation of the $RS_a$, that there is
a combinatorial method for calculating the action of $RS_a$ on a Schur function
when $m < \la_1$.  Let $ht_k(\mu)$ be the integer $i$ such that
$\mu \snake_k = (\mu_2 - 1, \mu_3 -1, \ldots, \mu_i-1,
\mu_1 + i - k, \mu_{i+1}, \ldots, \mu_{l(\mu)})$ is chosen 
to be a partition.  This amounts to removing the first $k$ cells from the border of
$\mu$.  If it is not possible to find such an $i$ such that
$\mu \snake_k$ is a partition then say that $\mu \snake_k$ is undefined.

\begin{cor}

Let $\nu = \la + (m+k)$ where $k \geq \la_1-m$ ($\nu$ is $\la$ resting on a sufficiently
long first row).

$$RS_a s_{\la} = (-1)^{ht_k(\nu)-1} s_{\nu \snake_k}$$

where it is assumed that $s_{\nu \snake_k} = 0$ if $\nu \snake_k$ does not exist.
\end{cor}
\vskip .3in

The proof of this corollary is not difficult, just a matter of showing
that the commutation relation of $RS_a RS_b$ agrees with this definition
of $\nu \snake_k$ and that the vanishing condition exists because $RS_a RS_{a+1} = 0$.
This definition and corollary are useful in showing that an expression for $(RS_a)^k$
can be reduced to a form that is very similar to the other vertex operators presented here.

\begin{lemma} \label{Smak}
$$(RS_a)^k = \sum_\la (-1)^{|\la|} s_{\la + a^k} s_{\la'}^\perp$$
with the convention that $s_{\la+a^k}$ is $0$ if $\la+a^k$ is undefined.
\end{lemma}
\vskip .3in

\noindent
{\bf Proof:}
By induction on $k$.  The statement agrees with Theorem \ref{bern} for $k=1$.

$$RS_a (RS_a)^k 
= \sum_{i \geq 0} (-1)^i h_{a+i} e_i^\perp \sum_{\la} (-1)^{|\la|} s_{\la+a^k} s_{\la'}^\perp$$
$e_i^\perp$ can be commuted with the Schur function to produce
$$= \sum_{i \geq 0} (-1)^i h_{a+i} \sum_{\la} (-1)^{|\la|} 
\sum_{j = 0}^i e_j^\perp( s_{\la+a^k}) e_{i-j}^\perp s_{\la'}^\perp$$

Interchange the order of all of the sums.
$$= \sum_{\la} \sum_{j \geq 0} \sum_{i \geq j} (-1)^{|\la|+i} h_{a+i}  
e_j^\perp( s_{\la+a^k}) e_{i-j}^\perp s_{\la'}^\perp$$

Make the substitution that $i \rightarrow i+j$, making the sum over all $i \geq 0$ and
expand the product $e_{i}^\perp s_{\la'}^\perp$.  The notation that
$\gamma \slash \la' \in \V_i$ means that $\gamma$ differs from $\la'$ by
a vertical $i$ strip ($\la'_j \geq \gamma_j \geq \la'_j+1$ and $|\gamma| = |\la|+i$). 
$$= \sum_{\la} \sum_{j \geq 0} \sum_{i \geq 0} (-1)^{|\la|+i+j} h_{a+i+j}  
e_j^\perp( s_{\la+a^k}) \sum_{\gamma \slash \la' \in \V_i} s_{\gamma}^\perp$$

Make the substitution $\gamma \rightarrow \gamma'$ so that the sum is over
all partitions $\gamma$ that differ from $\la$ by a horizontal $i$ strip
and rearrange the sums.
$$= \sum_{\la} \sum_{i \geq 0} \sum_{\gamma \slash \la \in \H_i} (-1)^{|\la|+i} 
\sum_{j \geq 0} (-1)^j h_{a+i+j} e_j^\perp( s_{\la+a^k}) s_{\gamma'}^\perp$$

Now it is only necessary to notice that the sum over $j$ is actually an application of
the Schur vertex operator acting exclusively on the symmetric function $s_{\la+a^k}$.
Switch the order of the sums over the partitions and expression becomes
$$= \sum_{\gamma } (-1)^{|\gamma|}
\sum_{i \geq 0} 
\sum_{\la: \gamma \slash \la \in \H_i}
RS_{a+i}( s_{\la+a^k}) s_{\gamma'}^\perp$$

There is a sign reversing involution on these terms so that only one term in the
sum over $i$ and $\la$ survives, namely, $s_{\gamma + a^{k+1}}$.
If $i=\gamma_1$ then $RS_{a+\gamma_1}(s_{\gamma-(\gamma_1)+a^k}) = s_{\gamma + a^{k+1}}$.

Take any partition $\la$ in this sum such that $\gamma \slash \lambda$ is a horizontal strip
of length less than $\gamma_1$.
If $RS_{a+i}(s_{\la+a^k}) = 0$, then this term does not contribute to the sum.
If $RS_{a+i}(s_{\la+a^k}) = s_{\nu + a^k}$ then $\nu = \la+(i+n) \snake_n$, where $n=\gamma_1-i$.
There is a combinatorial statement that can be made about partitions that satisfy this
condition, this is a lemma stated in \cite{Z} (Lemma 3.15, p.34).

\begin{lemma}
There exists an involution $I_\gamma^n$ on partitions $\mu$ such
that $\mu \slash \gamma$ is a horizontal $n$ strip, $\mu \snake_n$ exists and
$\gamma \neq \mu \snake_n$ with the property that $ht_n(I_\gamma^n(\mu)) =
ht_n(\mu) \pm 1$ and $\mu \snake_n = I_\gamma^n(\mu) \snake_n$.
\end{lemma}

This is exactly the situation here.  Set $\mu = \la + (i+n)$ then
$\mu \slash \gamma$ is a horizontal strip of size $|\mu| - |\gamma| = |\la|+i+n - |\gamma|
= n$.  The result then is that all terms cancel EXCEPT for the terms
such that $\gamma = \la+(i+n) \snake_n$ or $i=\gamma_1$ and 
$RS_{a+i}(s_{\la+a^k}) = s_{\gamma + a^{k+1}}$.

The sum therefore reduces to
$$=\sum_\gamma (-1)^{|\gamma|}
s_{\gamma + a^{k+1}} s_{\gamma'}^\perp$$
\endofproof
\vskip .3in

With this expression for the Schur function vertex operator, it is possible to reduce
the expression for the 'everything operator' that adds a column to the Schur functions
but is zero when the length of the indexing partition is larger than the height of
the column being added.

\begin{thm}
For $a,k \geq 0$, let $CS_{a^k} = \sum_\la (-1)^{|\la|} (RS_a)^k(s_\la) s_{\la'}^\perp$.
This operator has the property that $CS_{a^k} s_\la = s_{\la + a^k}$ if $l(\la) \leq k$ and 
$CS_{a^k} s_\la = 0$ for $l(\la) > 0$.
\end{thm}
\vskip .3in

\noindent
{\bf Proof:}

Take the expression for the everything operator that adds $a$ columns of height $k$
using the convention that $s_{\la+a^k}$ is zero whenever $\la+a^k$ is undefined.
$$E_{s_\la}^{s_{\la+a^k}} = \sum_\la s_{\la+a^k} \sum_\mu (-1)^{|\mu|} s_\mu s_{\mu'}^\perp 
s_{\la'}^\perp$$

The coefficients of the expansion of $s_{\mu} s_\la$ in terms of Schur functions
are well studied and
there exists formulas and combinatorial interpretations for their calculation.  The only
properties that we require here is that the coefficients in the
the expression $s_\mu s_\la = \sum_\nu c_{\la\mu}^\nu s_\nu$ have the property that
$c_{\la\mu}^\nu = c_{\la'\mu'}^{\nu'}$ and $s_\la^\perp s_\nu = \sum_\nu c_{\mu\la}^\nu s_\mu$.

$$= \sum_\la s_{\la+a^k} \sum_\mu (-1)^{|\mu|} s_\mu \sum_\nu c_{\la\mu'}^{\nu'} s_{\nu'}^\perp$$

Next, we rearrange the sums and make the substitution $c_{\la\mu'}^{\nu'} = c_{\la'\mu}^{\nu}$.
$$= \sum_\la s_{\la+a^k} \sum_\nu (-1)^{|\nu| - |\la|} \sum_\mu c_{\la'\mu}^{\nu} s_\mu
s_{\nu'}^\perp$$

Therefore the sum over $\mu$ is just an application of $s_{\la'}^\perp$ on $(s_{\nu})$ and the
sums can be rearranged.
$$= \sum_\nu (-1)^{|\nu|} \sum_\la (-1)^{|\la|} s_{\la+a^k} s_{\la'}^\perp (s_{\nu})
s_{\nu'}^\perp$$

The sum over $\la$ is now exactly an application of Lemma \ref{Smak}.
$$= \sum_\nu (-1)^{|\nu|} (RS_a)^k(s_{\nu}) s_{\nu'}^\perp$$
This is the expression given in the statement of the theorem.
\endofproof
\vskip .3in

The last of the 'eerie coincidences' of this article is that the $CS_{a^k}$ and $(RS_a)^k$
are related by a pair of formulas similar to the case of formulas (\ref{HErel}),
(\ref{EHrel}) and (\ref{CMrel}), (\ref{MCrel}).

\begin{equation}
CS_{a^k} = \sum_\la (-1)^{|\la|} (RS_a)^k(s_{\la}) s_{\la'}^\perp
\end{equation}
\begin{equation}
(RS_a)^k = \sum_\la (-1)^{|\la|} CS_{a^k}(s_{\la}) s_{\la'}^\perp
\end{equation}

\section{An application: the tableaux of bounded height}

One observation about the operator $CS_{a^k}$ that could have an interesting application is that
$CS_{0^k} s_\la = 0$ if $l(\la) > k$ and $CS_{0^k} s_\la = s_\la$ if $l(\la) \leq k$.
Knowing this and the commutation relation between $RS_a$ and $h_k$ allows us to calculate
the number of pairs of standard tableaux of the same shape of bounded height \cite{B}
$\sum_{\ontop{\la \vdash n}{l(\la)\leq
k}} f_\la^2$ (where $f_\la$ is the number of standard tableaux of shape $\la$).

\begin{prop}
Let $CP(n,k)$ be the collection of sequences of non-negative integers of
length $k$ such that the sum is $n$.
$$\sum_{\ontop{\la \vdash n}{l(\la)\leq
k}} f_\la^2 = \sum_{s \in CP(n,k)} \left( \begin{array}{c} n \\ s  \end{array} \right)
\frac{\prod_{i<j} (s_j+j - (s_i+i)) }{\prod_{i=1}^k (s_i + i -1)!} n!$$
\end{prop}

The formula follows by applying $CS_{0^k}$ to the symmetric function $h_1^n$ to arrive
at a formula for the symmetric function $\sum_{\ontop{\la \vdash n}{l(\la)\leq k+1}}
f_\la s_\la$.

\begin{lemma}
$$ CS_{0^k}(h_1^n)
= \sum_{\ontop{\la \vdash n}{l(\la) \leq k}} f_\la s_\la
= \sum_{s \in CP(n,k)} \left( \begin{array}{c} n\\ s\end{array} \right)
det \left| h_{s_j-j+i} \right|_{1 \leq i,j \leq k}$$
\end{lemma}
\vskip .3in

\noindent
{\bf Proof:}
Use the relation $RS_a h_k = h_k RS_a - h_{k-1} RS_{a+1}$,
$RS_a 1 = h_a$ and induction to calculate that

\begin{equation} \label{rsform}
RS_0^k (h_1^n) = \sum_{l=0}^n \sum_{s \in CP(n-l,k)} (-1)^{n-l} h_1^l \left(
\begin{array}{c} n \\ l,s \end{array} \right) det \left| h_{s_j-j+i} \right|_{1 \leq i,j \leq k}
\end{equation}

Using the relation that $s_\la^\perp (h_1^n) = \left( \begin{array}{c} n \\ |\la| \end{array}
\right) f_\la h_1^{n-|\la|}$ we have that
\begin{align*}
CS_{0^k}(h_1^n) &= \sum_\la (-1)^{|\la|} (RS_0)^k(s_\la) s_{\la'}^\perp (h_1^n) \\
&= \sum_\la (-1)^{|\la|} (RS_0)^k(s_\la) \left( \begin{array}{c} n \\ |\la| \end{array} \right)
f_\la h_1^{n-|\la|} \\
&= \sum_{i=0}^n \sum_{\la \vdash i} (-1)^{i} \left( \begin{array}{c} n \\ i \end{array}
\right) (RS_0)^k(f_\la s_\la) h_1^{n-i} \\
&= \sum_{i=0}^n (-1)^{i} \left( \begin{array}{c} n \\ i \end{array}
\right) (RS_0)^k(h_1^i) h_1^{n-i}
\end{align*}

Now using (\ref{rsform}) we can reduce this further to

$$=\sum_{m=0}^n \sum_{l=0}^m \sum_{s \in CP(m-l,k)} (-1)^{l}
\left( \begin{array}{c} n \\ m \end{array} \right) \left(
\begin{array}{c} m \\ l,s \end{array} \right) h_1^{n+l-m} 
det \left| h_{s_j-j+i} \right|_{1 \leq i,j \leq k}$$

Now switch the sums indexed by $l$ and $m$ and then make the replacement $m \rightarrow m+l$

$$=\sum_{l=0}^n \sum_{m=0}^{n-l} \sum_{s \in CP(m,k)} (-1)^{l}
\left( \begin{array}{c} n \\ m+l \end{array} \right) \left(
\begin{array}{c} m+l \\ l,s \end{array} \right)
h_1^{n-m} det \left| h_{s_j-j+i} \right|_{1 \leq i,j \leq k}$$
Now switch the sums back and rearrange the binomial coefficients
$$=\sum_{m=0}^n \sum_{l=0}^{n-m} \sum_{s \in CP(m,k)} (-1)^{l}
\left( \begin{array}{c} n \\ n-m,s \end{array} \right)
\left( \begin{array}{c} n-m \\ l \end{array} \right)
h_1^{n-m} det \left| h_{s_j-j+i} \right|_{1
\leq i,j \leq k}$$

Now the sum $\sum_{l=0}^{n-m} (-1)^{l} \left( \begin{array}{c} n-m \\ l \end{array} \right)$
will always be zero unless $n-m=0$ and if $n=m$ then it is $1$
and so the entire sum collapses to
$$=\sum_{s \in CP(n,k)} \left( \begin{array}{c} n \\ s \end{array} \right)
det \left| h_{s_j-j+i} \right|_{1 \leq i,j \leq k}$$ \endofproof
\vskip .3in

\noindent{\bf Proof:} (of proposition)

The proposition follows from this lemma with a little manipulation.  There is a linear
and multiplicative homomorphism that sends the symmetric functions to the space of
polynomials in one variable due to Gessel 
defined by $ \theta( h_n ) = x^n/n!$.  This homomorphism has the property
 that $\theta( s_\la ) = f_\la x^{|\la|}/|\la|!$.  The image of the formula 
in the lemma is then
$$\theta( CS_{0^k}(h_1^n) ) = \theta( \sum_{\ontop{\la \vdash n}{l(\la) \leq k}} f_\la s_\la )
= \sum_{\ontop{\la \vdash n}{l(\la) \leq k}} f_\la^2 \frac{x^n}{n!}$$

Therefore if we set $(a)_0 = 1$ and $(a)_i = a (a-1) \cdots (a-i+1)$ then we have 
(by making a slight transformation that reverses order of the sequence first...$j \rightarrow
n+1-j$, $i \rightarrow n+1-i$ and $s_i \rightarrow s_{n+1-i}$) that

\begin{equation*}
\sum_{\ontop{\la \vdash n}{l(\la)\leq
k}} f_\la^2 = \sum_{s \in CP(n,k)} \left( \begin{array}{c} n\\ s\end{array} \right)
det \left| \frac{(s_j+j-1)_{i-1}}{(s_j+j-1)!} \right|_{1 \leq i,j \leq k} n!
\end{equation*}

\begin{equation*}
\sum_{\ontop{\la \vdash n}{l(\la)\leq
k}} f_\la^2 = \sum_{s \in CP(n,k)} \left( \begin{array}{c} n\\ s\end{array} \right) \prod_{i=1}^k
\frac{1}{(s_j+j-1)!}
det \left| {(s_j+j-1)_{i-1}} \right|_{1 \leq i,j \leq k} n!
\end{equation*}

The determinant is a specialization of the Vandermonde determinant in the variables
$s_j + j -1$ so the formula reduces to the expression stated in the proposition. \endofproof

We note that in the case that $k=1$ this sum reduces to $1$ and in the case that $k=2$
we have that 
\begin{equation*}
\sum_{\ontop{\la \vdash n}{l(\la)\leq
2}} f_\la^2 = \sum_{j=0}^n \left( \begin{array}{c} n\\ j\end{array} \right)
\frac{n-2j+1}{(j)!(n-j+1)!} n!
= \sum_{j=0}^n \left( \begin{array}{c} n\\ j\end{array} \right)^2
\frac{n-2j+1}{ n-j+1 }
\end{equation*}
And this is an expression for the Catalan numbers.
It would be interesting to see if these expressions and equations
could be $q$ or $q,t$ anlogued.


\begin{thebibliography}{99}
\bibitem{B} F. Bergeron, L. Favreau, D. Krob, Conjectures on the enumeration
of tableaux of bounded height, {\it Discrete Mathematics}, 139 (1995) pp. 463-468

\bibitem{J} N. Jing,  Vertex operators and Hall-Littlewood symmetric functions,
{\it Adv. Math.}, {\bf 87} (1991), 226-248.

\bibitem{M} I.G. Macdonald,  "Symmetric Functions and Hall Polynomials,"
Oxford Mathematical Monographs,  Oxford UP, second edition, 1995.

\bibitem{Z} M. A. Zabrocki, A Macdonald Vertex Operator and
Standard Tableaux Statistics for the Two-Column $(q,t)$-Kostka Coefficients, 
{\it Electron. J. Combinat. 5}, R45 (1998), 46pp.

\bibitem{Ze} A. V. Zelevinsky. Representations of finite classical groups:
a Hopf algebra approach.  Springer Lecture Notes, 869 (1981)

\end{thebibliography}
\end{document}